\documentclass[a4paper, 10pt]{extarticle}
\usepackage{amsmath}
\usepackage{amssymb}
\usepackage{amsfonts}
\usepackage{graphicx}
\usepackage{subfigure}
\usepackage{authblk}
\usepackage{url}
\usepackage[textwidth=16cm,textheight=20cm]{geometry}
\newtheorem{theorem}{Theorem}

\newtheorem{example}[theorem]{Example}

\newtheorem{Remark}[theorem]{Remark}

\title{Singularities in Euler flows: multivalued solutions, shock waves, and phase transitions}
\author[1,\thanks{\textit{E-mail: }\texttt{valentin.lychagin@uit.no}}]{Valentin Lychagin}
\author[1,2\thanks{\textit{E-mail: }\texttt{mihail\underline{ }roop@mail.ru}}]{Mikhail Roop}
\affil[1]{V.A. Trapeznikov Institute of Control Sciences, Russian Academy of Sciences, 65 Profsoyuznaya Str., 117997 Moscow, Russia}
\affil[2]{Faculty of Physics, Lomonosov Moscow State University, Leninskie Gory 1, 119991 Moscow, Russia}

\begin{document}
\maketitle

\abstract{
In this paper, we analyze various types of critical phenomena in one-dimensional gas flows described by Euler equations. We give a geometrical interpretation of thermodynamics with a special emphasis on phase transitions. We use ideas from the geometrical theory of PDEs, in particular, symmetries and differential constraints to find solutions to the Euler system. Solutions obtained are multivalued, have singularities of projection to the plane of independent variables. We analyze the propagation of the shock wave front along with phase transitions.
}

\section{Introduction}
Various types of critical phenomena, such as singularities, discontinuities, wave fronts and phase transitions, have always been of interest from both mathematical \cite{Arn1,Arn2,Arn3} and practical \cite{Zeld} viewpoints. In the context of gases, discontinuous solutions to the Euler system, describing their motion, are usually treated as \textit{shock waves}. In the past decades such phenomena have widely been studied, see for instance, \cite{Huang} for the case of Chaplygin gases, \cite{Ros-Tab,Chat}, where the weak shocks are considered. It is worth mentioning also \cite{Polud1,Polud2}, where the influence of turbulence on shocks and detonations is emphasized.

This paper can be seen as a natural continuation of \cite{LR-ljm-shock}, where we have considered the case of ideal gas flows. Here, we elaborate the case of more complicated and at the same time more interesting from the singularity theory viewpoint model, the van der Waals model, known as one of the most popular in the description of phase transitions. So, singularities of shock wave type that can be viewed as in some sense singular solutions to the Euler system are analyzed together with singularities of purely thermodynamic nature, phase transitions. Our approach to finding and investigating such phenomena is essentially based on the geometric theory of PDEs \cite{KLR,KVL,KrVin,Olver,Ovs}. Namely, we find a class of multivalued solutions to the Euler system (see also \cite{Tun}), and singularities of their projection to the plane of independent variables are exactly what drives the appearance of the shock wave \cite{LychSing}. Similar ideas were used in a series of works \cite{AKL-dan,AKL-ifac,AKL-gsa}, where multivalued solutions to filtration equations were obtained along with analysis of shocks. To find such solutions we use the idea of adding a differential constraint to the original PDE in such a way that the resulting overdetermined system of PDEs is compatible \cite{KrugLych}. The same concepts were also used in \cite{Eiv}, where a general solution to the Hunter-Saxton equation was found, in \cite{Lych-Yum-1}, where the two-dimensional Euler system was considered, and in \cite{Lych-Yum-2}, where this approach was applied to the Khokhlov-Zabolotskaya equation.

The paper is organized as follows. Section \ref{sec1} is preliminary, we describe there the necessary concepts from thermodynamics. In Section \ref{sec3}, we analyze a multivalued solution to Euler equations and its singularities, including shock waves and phase transitions. In the last section, we discuss the results. The essential computations for this paper were made with the DifferentialGeometry package \cite{And} in Maple.
\section{Thermodynamics}
\label{sec1}
In this section, we give necessary concepts from thermodynamics. As we shall see, geometrical interpretation of thermodynamic states allows one to use Arnold's ideas from the theory of Legendrian and Lagrangian singularities \cite{Arn1,Arn2,Arn3}, which are crucial in description of phase transitions. The geometrical approach to thermodynamics was initiated already by Gibbs \cite{Gibbs}, and was further developed in, for example, \cite{Mrug,Rup}, and recently in \cite{Lych}. For more detailed analysis we also refer to \cite{LR-ljm}.
\subsection{Legendrian and Lagrangian manifolds}
Consider the contact space $\left(\mathbb{R}^{5},\theta\right)$ with coordinates $(s,e,\rho,p,T)$ standing for specific entropy, specific inner energy, density, pressure and temperature. The contact structure $\theta$ is given by

\begin{equation}
\label{cont}
\theta=T^{-1}de-ds-pT^{-1}\rho^{-2}d\rho.
\end{equation}

Then, a \textit{thermodynamic state} is a Legendrian manifold $\widehat{L}\subset\left(\mathbb{R}^{5},\theta\right)$, i.e. $\theta\left|_{\widehat{L}}\right.=0$ and $\dim\widehat{L}=2$. From the physical viewpoint this means that the first law of thermodynamics holds on $\widehat{L}$. Due to (\ref{cont}), it is natural to choose $(e,\rho)$ as coordinates on $\widehat{L}$. Then, a two-dimensional manifold $\widehat{L}\subset\left(\mathbb{R}^{5},\theta\right)$ is given by

\begin{equation}
\label{leg}
\widehat{L}=\left\{s=S(e,\rho),\, T=\frac{1}{S_{e}},\, p=-\rho^{2}\frac{S_{\rho}}{S_{e}}\right\},
\end{equation}

where the function $S(e,\rho)$ specifies the dependence of the specific entropy on $e$ and $\rho$.

Note that determining a Legendrian manifold $\widehat{L}$ by means of (\ref{leg}) requires the knowledge of $S(e,\rho)$, while in experiments one usually gets relations between pressure, density and temperature. To this reason, we get rid of the specific entropy $s$ by means of projection $\pi\colon\mathbb{R}^{5}\to\mathbb{R}^{4}$, $\pi(s,e,\rho,p,T)=(e,\rho,p,T)$ and consider an immersed Lagrangian manifold $\pi\left(\widehat{L}\right)=L\subset\left(\mathbb{R}^{4},\Omega\right)$ in a symplectic space $\left(\mathbb{R}^{4},\Omega\right)$, where the structure symplectic form $\Omega$ is

\begin{equation*}
\Omega=d\theta=d(T^{-1})\wedge de-d(pT^{-1}\rho^{-2})\wedge d\rho.
\end{equation*}

Then, one can treat thermodynamic state manifolds as Lagrangian manifolds $L\subset\left(\mathbb{R}^{4},\Omega\right)$, i.e. $\Omega\left|_{L}\right.=0$. In coordinates $(T,\rho)$, a thermodynamic Lagrangian manifold $L$ is given by two functions

\begin{equation}
\label{lag}
L=\left\{p=P(T,\rho),\,e=E(T,\rho)\right\}.
\end{equation}

Since $\Omega\left|_{L}\right.=0$, the functions $P(T,\rho)$ and $E(T,\rho)$ are not arbitrary, but are related by

\begin{equation}
\label{Pbr}
\left[p-P(T,\rho),e-E(T,\rho)\right]\left|_{L}\right.=0,
\end{equation}

where $[f,g]$ is a Poisson bracket of functions $f$ and $g$ on $\left(\mathbb{R}^{4},\Omega\right)$ uniquely defined by the relation

\begin{equation*}
[f,g]\,\Omega\wedge\Omega=df\wedge dg\wedge\Omega.
\end{equation*}

Equation (\ref{Pbr}) forces the following relation between $P(T,\rho)$ and $E(T,\rho)$: $(-\rho^{-2}T^{-1}P)_{T}=(T^{-2}E)_{\rho}$, and therefore the following theorem is valid:
\begin{theorem}
The Lagrangian manifold $L$ is given by means of \textit{the Massieu-Planck potential} $\phi(\rho,T)$

\begin{equation}
\label{LagPhi}
p=-\rho^{2}T\phi_{\rho},\quad e=T^{2}\phi_{T}.
\end{equation}

\end{theorem}
\begin{Remark}
Having given the Lagrangian manifold $L$ by means of (\ref{lag}), one can find out the entropy function $S(e,\rho)$ solving the overdetermined system

\begin{equation*}
T=\frac{1}{S_{e}},\, p=-\rho^{2}\frac{S_{\rho}}{S_{e}}
\end{equation*}

with compatibility condition (\ref{Pbr}).
\end{Remark}

\subsection{Riemannian structures, singularities, phase transitions}
There is one more important structure arising, as it was shown in \cite{Lych}, from measurement approach to thermodynamics. Indeed, if one considers equilibrium thermodynamics as a theory of measurement of random vectors, whose components are inner energy and volume $v=\rho^{-1}$, one drives to the universal quadratic form on $(\mathbb{R}^{4},\Omega)$ of signature $(2,2)$:

\begin{equation*}
\kappa=d(T^{-1})\cdot de-\rho^{-2}d(pT^{-1})\cdot d\rho,
\end{equation*}

where $\cdot$ is the symmetric product, and areas on $L$, where the restriction $\kappa|_{L}$ of $\kappa$ to $L$ is negative, are those where the variance of a random vector $(e,v=\rho^{-1})$ is positive \cite{Lych,KLR-ent}. Using (\ref{LagPhi}), we get

\begin{equation}
\label{kappaPhi}
\kappa|_{L}=-(2T^{-1}\phi_{T}+\phi_{TT})dT\cdot dT+(2\rho^{-1}\phi_{\rho}+\phi_{\rho\rho})d\rho\cdot d\rho,
\end{equation}

and taking into account (\ref{LagPhi}), we conclude that the condition of positive variance is satisfied at points on $L$, where

\begin{equation*}
e_{T}>0,\quad p_{\rho}>0,
\end{equation*}

which is known as the condition of the thermodynamic stability.

Let us now explore singularities of Lagrangian manifolds. We will be interested in the singularities of their projection to the plane of intensive variables $(p,T)$, i.e. points where the form $dp\wedge dT$ degenerates. We will assume that extensive variables $(e,\rho)$ may serve as global coordinates on $L$, i.e. the form $de\wedge d\rho$ is non-degenerated everywhere. The set where $dp\wedge dT=0$ coincides with that where $2\rho^{-1}\phi_{\rho}+\phi_{\rho\rho}=0$, or, equivalently, where the from $\kappa|_{L}$ degenerates. A manifold $L$ turns out to be divided into submanifolds $L_{i}$, where both $(e,\rho)$ and $(p,T)$ may serve as coordinates, or, equivalently, the form (\ref{kappaPhi}) is non-degenerated. Such $L_{i}$ are called \textit{phases}. Additionally, those of $L_{i}$, where (\ref{kappaPhi}) is negative, are called \textit{applicable phases}. Thus we end up with the observation that singularities of projection of thermodynamic Lagrangian manifolds are related with the theory of phase transitions. Indeed, by a \textit{phase transition} of the first order we mean a jump from one applicable state to another, governed by the conservation of intensive variables $p$, $T$, and specific Gibbs potential

\begin{equation*}
\gamma=e-Ts+p/\rho,
\end{equation*}

which in terms of the Massieu-Planck potential is expressed as $\gamma=-T(\phi+\rho\phi_{\rho})$ \cite{LR-ljm}. Consequently, to find the points of phase transition, one needs to solve the system

\begin{equation}
\label{ph-tr-cond}
p=-\rho_{1}^{2}T\phi_{\rho}(T,\rho_{1}),\quad p=-\rho_{2}^{2}T\phi_{\rho}(T,\rho_{2}),\quad \phi(T,\rho_{1})+\rho_{1}\phi_{\rho}(T,\rho_{1})=\phi(T,\rho_{2})+\rho_{2}\phi_{\rho}(T,\rho_{2}),
\end{equation}

where $p$ and $T$ are the pressure and the temperature of the phase transition, and $\rho_{1}$ and $\rho_{2}$ are the densities of gas and liquid phases.
\begin{example}[Ideal gas]
The simplest example of a gas is an ideal gas model. In this case the Legendrian manifold is given by

\begin{equation}
\label{leg-ideal}
\widehat{L}=\left\{p=R\rho T,\, e=\frac{n}{2}RT,\, s=R\ln\left(\frac{T^{n/2}}{\rho}\right)\right\},
\end{equation}

where $R$ is the universal gas constant, and $n$ is the degree of freedom. The differential quadratic form $\kappa|_{L}$ is

\begin{equation*}
\kappa|_{L}=-\frac{Rn}{2}\frac{dT^{2}}{T^{2}}-R\rho^{-2}d\rho^{2}.
\end{equation*}

It is negative definite on the entire $\widehat{L}$, and there are no phase transitions, nor are there singularities of projection of $\widehat{L}$ to the $p-T$ plane.
\end{example}
\begin{example}[van der Waals gas]
To define the Legendrian manifold for van der Waals gases we will use reduced state equations:

\begin{equation}
\label{leg-vdW}
\widehat{L}=\left\{p=\frac{8T\rho}{3-\rho}-3\rho^{2},\, e=\frac{4nT}{3}-3\rho,\,s=\ln\left(T^{4n/3}(3\rho^{-1}-1)^{8/3}\right)\right\}.
\end{equation}

The differential quadratic form $\kappa|_{L}$ is

\begin{equation*}
\kappa|_{L}=-\frac{4n}{3T^{2}}dT^{2}+\frac{6(\rho^{3}-6\rho^{2}-4T+9\rho)}{\rho^{2}T(\rho-3)^{2}}d\rho^{2}.
\end{equation*}

In this case it changes its sign, the manifold $\widehat{L}$ has a singularity of cusp type. The singular set of $\widehat{L}$, called also \textit{caustic}, and the curve of phase transition are shown in Figure 1.
\begin{figure}[ht!]
\centering
\subfigure[]{\includegraphics[width=0.4\linewidth]{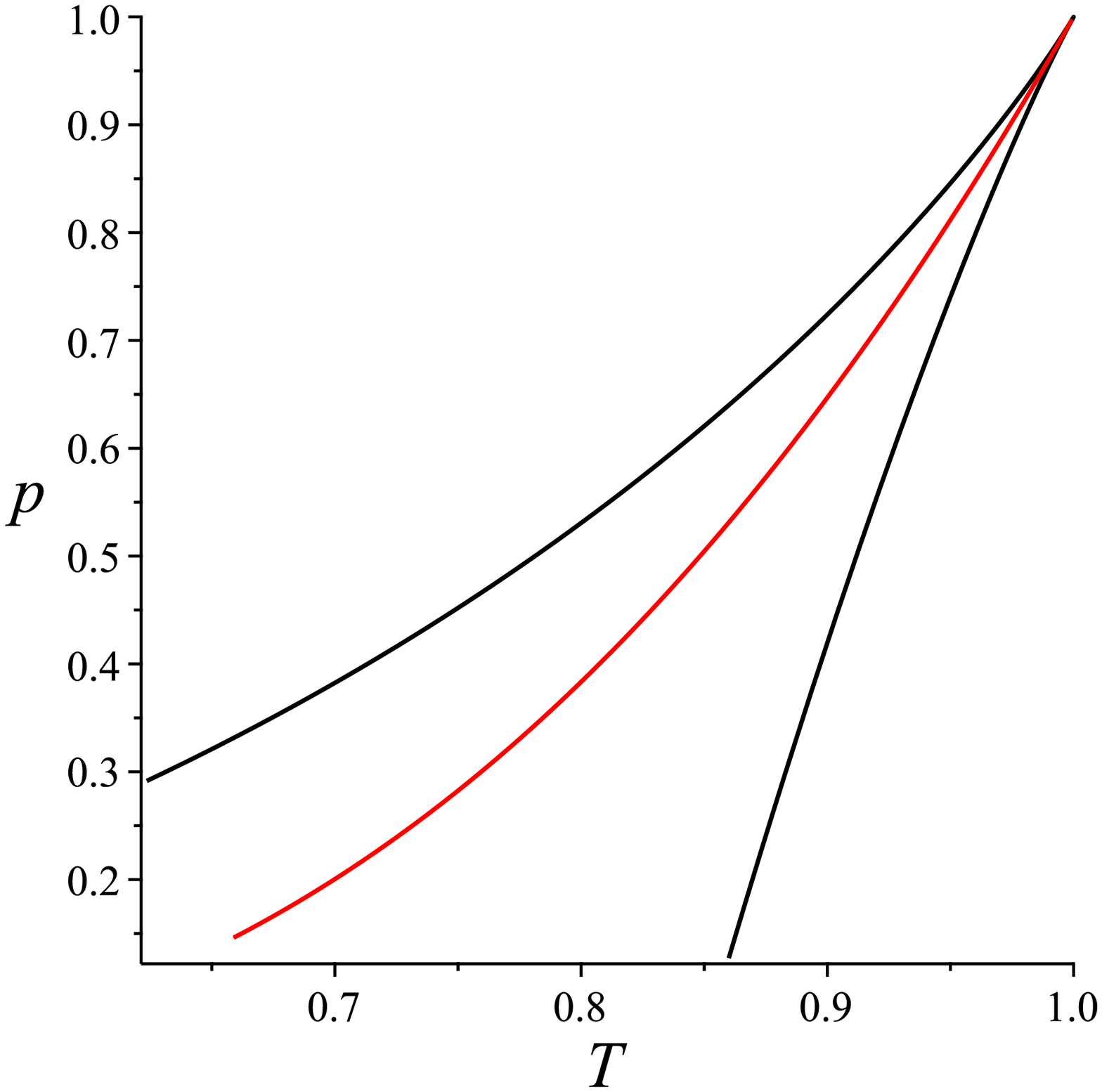} \label{lam1} }
\hspace{4ex}
\subfigure[]{ \includegraphics[width=0.5\linewidth]{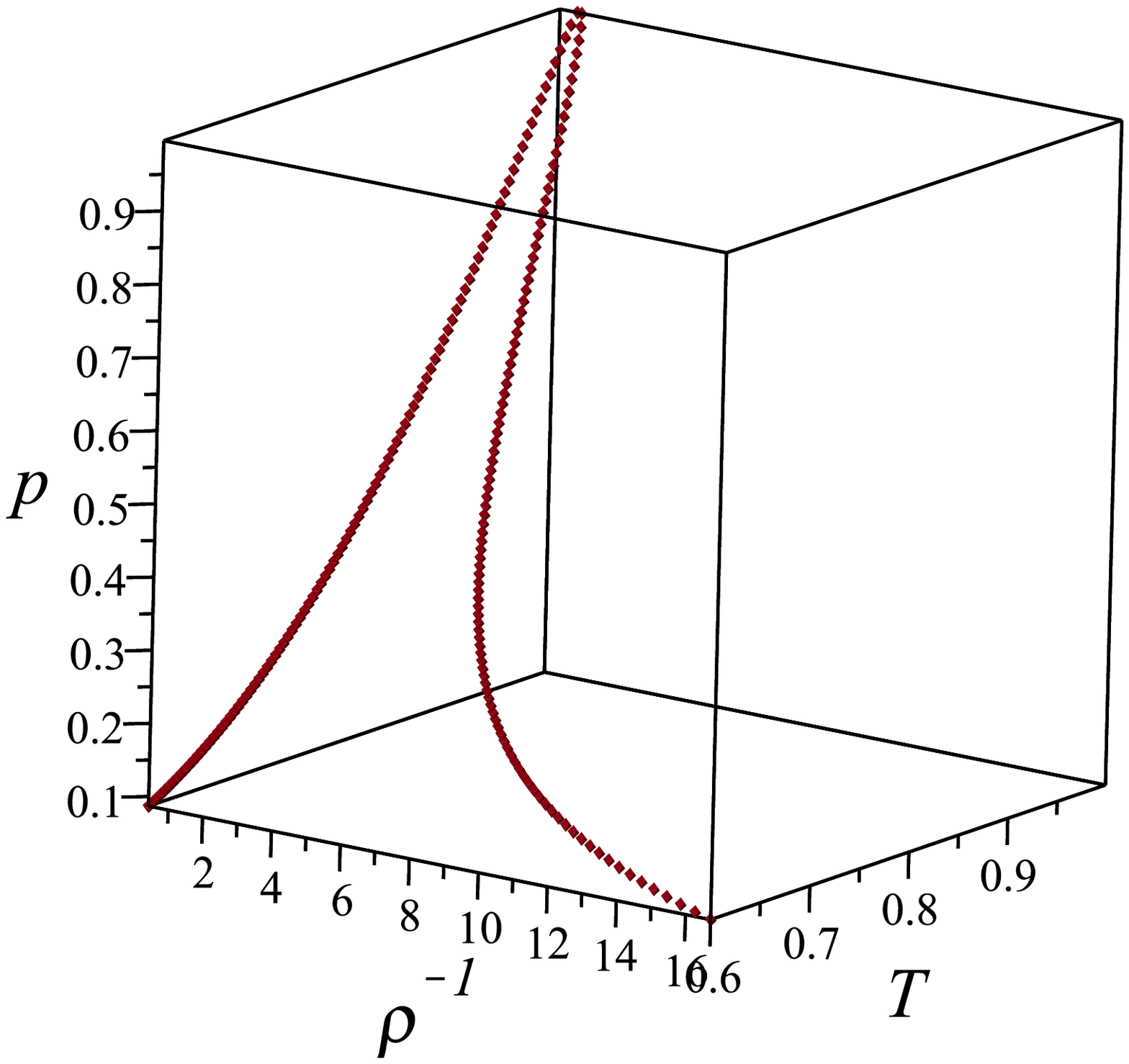} \label{lam2} }
\hspace{4ex}
\caption{\footnotesize{Singularities of the van der Waals Legendrian manifold. Caustic (black line) and phase transition curve (red line) in coordinates $(p,T)$ \subref{lam1}; the curve of phase transition in $(p,\rho,T)$ \subref{lam2}. Gas phase on the right of the curve, liquid phase on the left of the curve, wet steam inside the curve.}}
\label{dens-mex}
\end{figure}

\end{example}
\section{Euler equations}
\label{sec3}
In this paper we will study non-stationary, one-dimensional flows of gases, described by the following system of differential equations:
\begin{itemize}
\item Conservation of momentum

\begin{equation}
\label{Euler}
\rho(u_{t}+uu_{x})=-p_{x},
\end{equation}

\item Conservation of mass

\begin{equation}
\label{mass}
\rho_{t}+(\rho u)_{x}=0,
\end{equation}

\item Conservation of entropy along the flow

\begin{equation}
\label{energy}
s_{t}+us_{x}=0.
\end{equation}

\end{itemize}
Here $u(t,x)$ is the flow velocity, $\rho(t,x)$ is the density of the medium, and $s(t,x)$ is the specific entropy. System (\ref{Euler})-(\ref{energy}) is incomplete. It becomes complete once extended by equations of thermodynamic state (\ref{leg}). We will be interested in \textit{homentropic flows}, i.e. those with $s(t,x)=s_{0}$. On the one hand, this assumption satisfies (\ref{energy}) identically, on the other hand, it allows us to express all the thermodynamic variables in terms of $\rho$. Indeed, the entropy $s$ has the following expression in terms of the Massieu-Planck potential $\phi(T,\rho)$: $s=\phi+T\phi_{T}$ \cite{LR-ljm}. Putting $s=s_{0}$, we get an equation $s_{0}=\phi+T\phi_{T}$, which determines $T(\rho)$ uniquely, since its derivative with respect to $T$ is positive due to the negativity of $\kappa|_{L}$. Substituting $T(\rho)$ into (\ref{lag}), one gets $p=p(\rho)$. Thus we end up with the following two-component system of PDEs:

\begin{equation}
\label{sysEu}
u_{t}+uu_{x}+A(\rho)\rho_{x}=0,\quad \rho_{t}+(\rho u)_{x}=0,
\end{equation}

where $A(\rho)=p^{\prime}(\rho)/\rho$.

We do not specify the function $A(\rho)$ yet, we shall do this while solving (\ref{sysEu}).

\subsection{Finding solutions}
To find solutions to system (\ref{sysEu}), we use the idea of adding a differential constraint to (\ref{sysEu}), compatible with the original system. It is worth mentioning that a solution is an integral manifold of the Cartan distribution on (\ref{sysEu}) (see \cite{KLR,KVL,KrVin} for details). This geometrical interpretation of a solution to a PDE allows to find ones in the form of manifolds, which, in general, may not be globally given by functions. This approach gives rise to investigation of singularities in a purely geometrical manner, which is shown in this paper.

 In general, finding differential constraints is not a trivial problem. But having found ones, the problem of finding solutions is reduced to the integration of a completely integrable Cartan distribution of the resulting compatible overdetermined system. In case the Cartan distribution has a solvable transversal symmetry algebra, whose dimension equals the codimension of the Cartan distribution, we are able to get explicit solutions in quadratures by applying the Lie-Bianchi theorem (for details we refer to \cite{KLR,KVL,KrVin}).

We will look for a differential constraint compatible with (\ref{sysEu}) in the form of a quasilinear equation

\begin{equation}
\label{constr}
u_{x}-\rho_{x}\left(\alpha(\rho)u+\beta(\rho)\right)=0,
\end{equation}

where functions $\alpha(\rho)$ and $\beta(\rho)$ are to be determined. We will denote system (\ref{sysEu})-(\ref{constr}) by $\mathcal{E}$.
\begin{theorem}
\label{thm1}
System (\ref{sysEu})-(\ref{constr}) is compatible if

\begin{equation}
\label{constr-solut}
\alpha(\rho)=-\frac{1}{\rho(C_{3}\rho-1)},\quad\beta(\rho)=\frac{C_{2}}{\rho(C_{3}\rho-1)},\quad A(\rho)=C_{1}+\frac{C_{5}}{\rho^{3}}\left(C_{3}+\frac{C_{7}}{\rho}\right)^{C_{6}},
\end{equation}

where $C_{i}$ are constants.

\end{theorem}
The proof of Theorem \ref{thm1} is more technical rather than conceptual. First of all, we lift system (\ref{sysEu})-(\ref{constr}) to the space of 3-jets $J^{3}(\mathbb{R}^{2})$ by applying total derivatives

\begin{eqnarray*}
D_{t}&=&\partial_{t}+u_{t}\partial_{u}+\rho_{t}\partial_{\rho}+u_{tt}\partial_{u_{t}}+\rho_{tt}\partial_{\rho_{t}}+\ldots,\\
D_{x}&=&\partial_{x}+u_{x}\partial_{u}+\rho_{x}\partial_{\rho}+u_{xx}\partial_{u_{x}}+\rho_{xx}\partial_{\rho_{x}}+\ldots.
\end{eqnarray*}

to equations of $\mathcal{E}$ the required number of times consequently. The resulting system $\mathcal{E}_{3}\subset J^{3}(\mathbb{R}^{2})$, consisting of equations only of the third order, contains 9 equations for 8 variables of purely third order, i.e. $u_{ttt}$, $u_{xxx}$, $u_{txx}$ and so on. Eliminating them from $\mathcal{E}_{3}$, we get 7 relations (6 obtained by lifting $\mathcal{E}$ to $J^{2}(\mathbb{R}^{2}$) plus one remained from eliminations of third-order variables). Again, we eliminate all the variables of the second order and we get four relations of the first order. Eliminating $u_{x}$, $u_{t}$ and $\rho_{t}$ we end up with the expression of the form $\rho_{x}^{3}G(\rho,u)=0$, where $G(\rho,u)$ is a polynomial in $u$, whose coefficients are ODEs on $\alpha(\rho)$, $\beta(\rho)$ and $A(\rho)$, solving which we get (\ref{constr-solut}). It is worth saying that these computations are algebraic and are well suited for computer algebra systems.
\begin{Remark}
Using (\ref{leg-ideal}) and (\ref{leg-vdW}), one can show that the function $A(\rho)=p^{\prime}(\rho)/\rho$ given in (\ref{constr-solut}) corresponds to that of
\begin{itemize}
\item ideal gas in case of

\begin{equation*}
C_{1}=C_{3}=0,\quad C_{5}=R\left(1+\frac{2}{n}\right)\exp\left(\frac{2s_{0}}{Rn}\right),\quad C_{6}=-2-\frac{2}{n},\quad C_{7}=1.
\end{equation*}

\item van der Waals gas in case of

\begin{equation}
\label{vdW-const}
C_{1}=-6,\quad C_{3}=-1,\quad C_{5}=24\left(1+\frac{2}{n}\right)\exp\left(\frac{3s_{0}}{4n}\right),\quad C_{6}=-2-\frac{2}{n},\quad C_{7}=3.
\end{equation}
\end{itemize}
The case of ideal gases was thoroughly investigated in \cite{LR-ljm-shock}. Here, we are interested in the case of van der Waals gases.
\end{Remark}

Summarizing, we have a compatible overdetermined system of PDEs

\begin{equation*}
\mathcal{E}=\left\{F_{1}=u_{t}+uu_{x}+A(\rho)\rho_{x}=0,\, F_{2}=\rho_{t}+(\rho u)_{x}=0,\, F_{3}=u_{x}-\rho_{x}\left(\alpha(\rho)u+\beta(\rho)\right)=0\right\}\subset J^{1}(\mathbb{R}^{2}),
\end{equation*}

where functions $\alpha(\rho)$, $\beta(\rho)$ and $A(\rho)$ are specified in (\ref{constr-solut}). This system is a smooth manifold $\mathcal{E}$ in the space of 1-jets $J^{1}(\mathbb{R}^{2})$ of functions on $\mathbb{R}^{2}$. Since $\dim J^{1}(\mathbb{R}^{2})=8$, and $\mathcal{E}$ consists of 3 relations on $J^{1}(\mathbb{R}^{2})$, $\dim\mathcal{E}=5$. The dimension of the Cartan distribution $\mathcal{C}_{\mathcal{E}}$ on $\mathcal{E}$ equals 2, therefore $\mathrm{codim}\,\mathcal{C}_{\mathcal{E}}=3$. Let us choose $(t,x,u,\rho,\rho_{x})$ as internal coordinates on $\mathcal{E}$. Then the Cartan distribution $\mathcal{C}_{\mathcal{E}}$ is generated by differential 1-forms

\begin{eqnarray}
\label{Cart-forms}
\omega_{1}&=&du-u_{x}dx-u_{t}dt,\\
\omega_{2}&=&d\rho-\rho_{x}dx-\rho_{t}dt,\\
\omega_{3}&=&d\rho_{x}-\rho_{xx}dx-\rho_{xt}dt,
\end{eqnarray}

where $\rho_{xx}$, $\rho_{xt}$, $u_{t}$, $u_{x}$, $\rho_{t}$ are expressed due to $\mathcal{E}$ and its prolongation \\$\mathcal{E}_{2}=\left\{D_{t}(F_{1})=0,\,D_{t}(F_{2})=0,\,D_{t}(F_{3})=0,\,D_{x}(F_{1})=0,\,D_{x}(F_{2})=0,\,D_{x}(F_{3})=0\right\}$:

\begin{equation}
\label{expr1}
\rho_{xx}=\frac{\rho_{x}^{2}\left(\rho(C_{3}\rho-1)^{3}A^{\prime}+(C_{3}\rho-1)^{2}A+3C_{3}(C_{2}-u)^{2}\right)}{(C_{3}\rho-1)\left((C_{2}-u)^{2}-A\rho(C_{3}\rho-1)^{2}\right)},\quad \rho_{t}=\frac{\rho_{x}(C_{3}\rho u+C_{2}-2u)}{1-C_{3}\rho},\end{equation}
\begin{equation}\label{expr2}u_{x}=\frac{\rho_{x}(C_{2}-u)}{\rho(C_{3}\rho-1)},\quad u_{t}=-\frac{\rho_{x}\left(A\rho(C_{3}\rho-1)+u(C_{2}-u)\right)}{\rho(C_{3}\rho-1)},\end{equation}

\begin{equation}
\label{expr3}
\begin{split}
\rho_{xt}=&\frac{\rho_{x}^{2}}{\rho(C_{3}\rho-1)^{2}\left(A\rho(C_{3}\rho-1)^{2}-(C_{2}-u)^{2}\right)}\left(\rho^{2}(C_{3}\rho-1)^{3}(C_{3}\rho u+C_{2}-2u)A^{\prime}\right.+{}\\&+\left.\rho A(C_{3}\rho-1)^{2}(C_{3}\rho u+3C_{2}-4u)+(C_{2}-u)^{2}(3C_{3}^{2}\rho^{2}u+3C_{3}\rho(C_{2}-2u)-2C_{2}+2u)\right),
\end{split}
\end{equation}

where $A(\rho)$ is given by (\ref{constr-solut}). We shall look for integrals of the distribution (\ref{Cart-forms})-(\ref{expr3}), which give us an (implicit) solution to (\ref{sysEu})-(\ref{constr}).
\begin{theorem}
The distribution (\ref{Cart-forms})-(\ref{expr3}) is a completely integrable distribution with a 3-dimensional Lie algebra $\mathfrak{g}$ of transversal infinitesimal symmetries generated by vector fields

\begin{equation*}
X_{1}=t\partial_{t}+x\partial_{x}-\rho_{x}\partial_{\rho_{x}},\quad X_{2}=\partial_{t},\quad X_{3}=\partial_{x}
\end{equation*}

with brackets $[X_{1},X_{3}]=-X_{3}$, $[X_{1},X_{2}]=-X_{2}$, $[X_{2},X_{3}]=0$.

The Lie algebra $\mathfrak{g}$ is solvable, and its sequence of derived algebras is

\begin{equation*}
\mathfrak{g}=\langle X_{1}, X_{2}, X_{3}\rangle\supset\langle X_{2}, X_{3}\rangle\supset 0.
\end{equation*}

\end{theorem}
Thus the Lie-Bianchi theorem \cite{KLR,KVL,KrVin} can be applied to integrate (\ref{Cart-forms})-(\ref{expr3}).

Let us choose another basis $\langle\varkappa_{1},\varkappa_{2},\varkappa_{3}\rangle$ in $\mathcal{C}_{\mathcal{E}}$ by the following way:

\begin{equation*}
\begin{pmatrix}
\varkappa_{1}\\\varkappa_{2}\\\varkappa_{3}
\end{pmatrix}
=
\begin{pmatrix}
\omega_{1}(X_{1}) & \omega_{1}(X_{2}) & \omega_{1}(X_{3})\\
\omega_{2}(X_{1}) & \omega_{2}(X_{2}) & \omega_{2}(X_{3})\\
\omega_{3}(X_{1}) & \omega_{3}(X_{2}) & \omega_{3}(X_{3})
\end{pmatrix}^{-1}
\begin{pmatrix}
\omega_{1}\\\omega_{2}\\\omega_{3}
\end{pmatrix}
.
\end{equation*}

Due to the structure of the symmetry Lie algebra $\mathfrak{g}$, the form $\varkappa_{1}$ is closed \cite{KLR,KrVin}, and therefore locally exact, i.e. $\varkappa_{1}=dQ_{1}$, where $Q_{1}\in C^{\infty}(J^{1})$, while restrictions $\varkappa_{2}|_{M_{1}}$ and $\varkappa_{3}|_{M_{1}}$ to the manifold $M_{1}=\left\{Q_{1}=\mathrm{const}\right\}$ are closed and locally exact too. Integrating the differential 1-form $\varkappa_{1}$ we observe that variables $u$, $\rho$, $t$, $x$ can be chosen as local coordinates on $M_{1}$ and

\begin{equation*}
M_{1}=\left\{\rho_{x}=\frac{\alpha_{1}\rho^{2}(C_{3}\rho-1)}{\rho A(C_{3}\rho-1)^{2}-(C_{2}-u)^{2}}\right\},
\end{equation*}

where $\alpha_{1}$ is a constant. Integrating restrictions $\varkappa_{2}|_{M_{1}}$ and $\varkappa_{3}|_{M_{1}}$, we get two more relations that give us a solution to (\ref{sysEu})-(\ref{constr}) implicitly:

\begin{equation}
\label{quad1}
t+\alpha_{2}+\frac{C_{2}-u}{\alpha_{1}\rho}+\frac{C_{3}u}{\alpha_{1}}=0,
\end{equation}

and

\begin{equation}
\label{quad2}
\begin{split}
0&=x+\alpha_{3}+\frac{1}{\alpha_{1}}\left(C_{1}\ln\rho-C_{1}C_{3}\rho+\frac{C_{3}u^{2}}{2}+\frac{u(C_{2}-u)}{\rho}-C_{5}\left(C_{3}+\frac{C_{7}}{\rho}\right)^{C_{6}+1}\cdot\right.{}\\&\cdot
\left.\frac{2\rho^{2}C_{3}^{2}-C_{7}^{2}(C_{6}+1)(C_{3}\rho(C_{6}+3)-C_{6}-2)+C_{3}C_{7}\rho(C_{3}\rho(C_{6}+3)-2C_{6}-2)}{(C_{6}+1)(C_{6}+2)(C_{6}+3)C_{7}^{3}\rho^{2}}\right),
\end{split}
\end{equation}

where we have already substituted $A(\rho)$ from (\ref{constr-solut}), and $\alpha_{2}$, $\alpha_{3}$ are constants. The graph of a multivalued solution for the density is shown in Figure 2. We used substitution (\ref{vdW-const}), where $C_{5}=240$, $n=3$, together with $C_{2}=1$, $\alpha_{1}=1$, $\alpha_{2}=2$, $\alpha_{3}=1$.

\begin{figure}[h!]

\begin{minipage}[h]{0.5\linewidth}
\center{\includegraphics[scale=0.3]{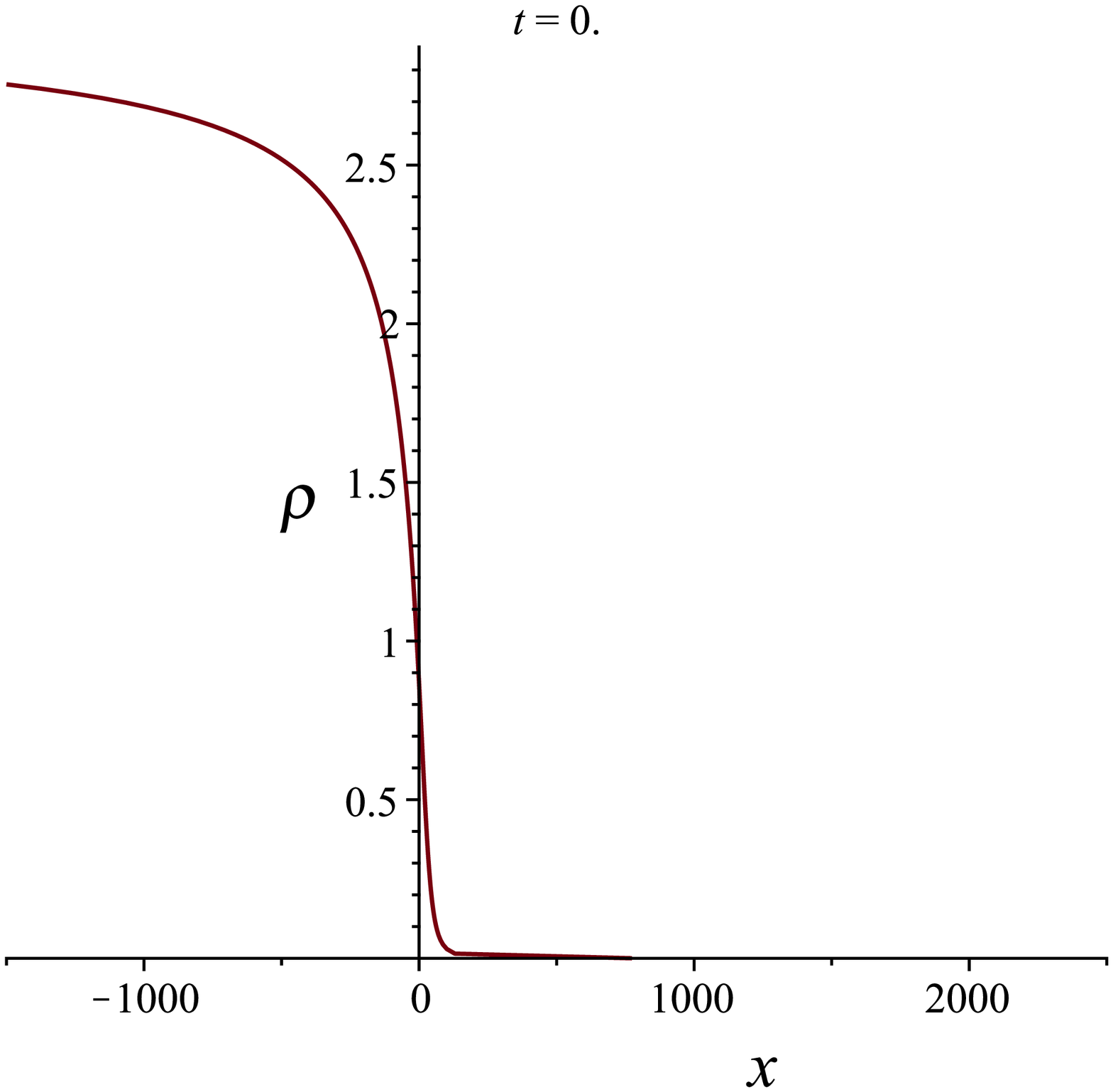}}
\end{minipage}
\hfill
\begin{minipage}[h]{0.5\linewidth}
\center{\includegraphics[scale=0.3]{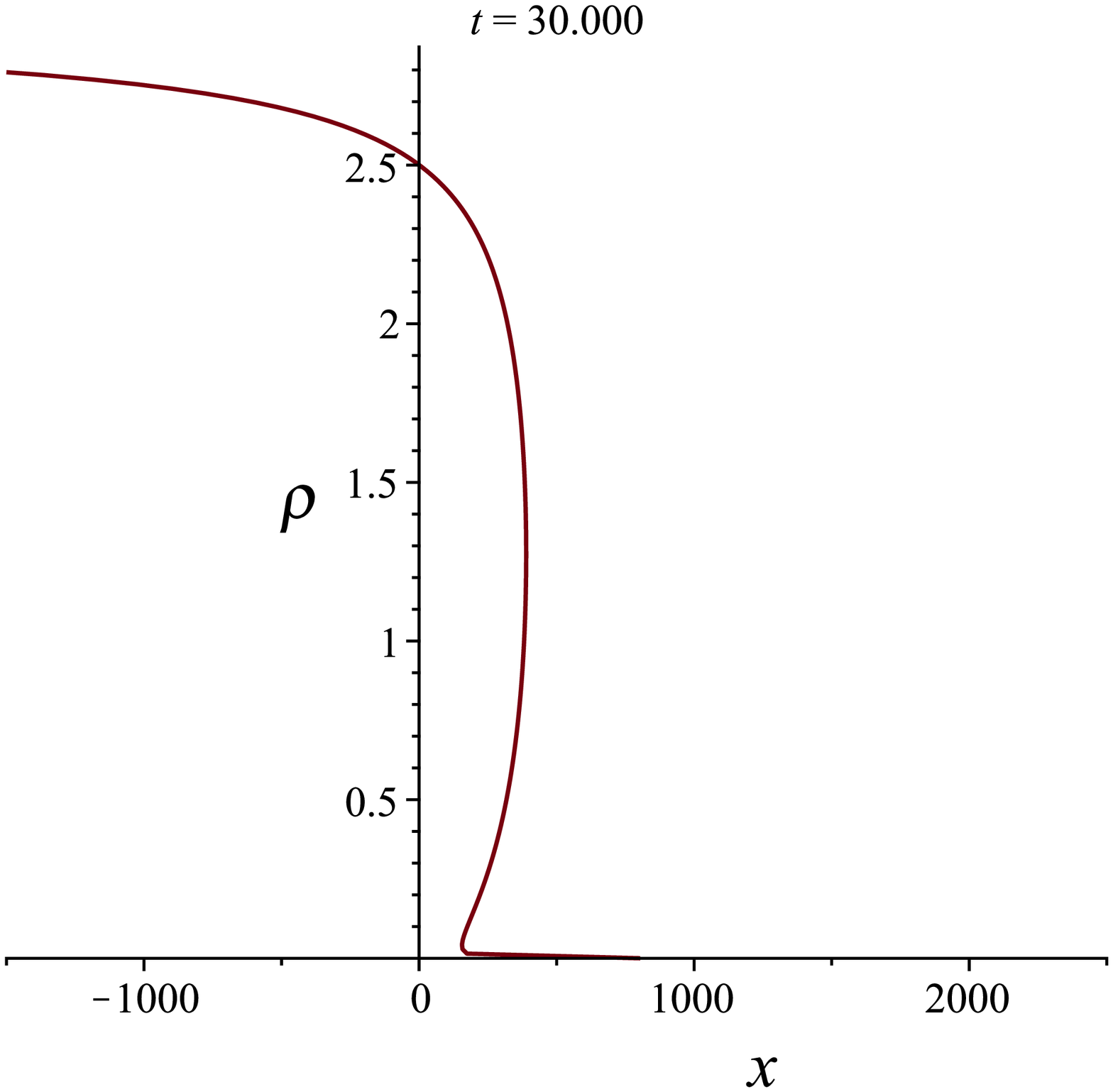}}
\end{minipage}
\hfill
\caption{Graph of the density in case of $n=3$ for time moments $t=0$, $t=30$}
\end{figure}

\subsection{Caustics and shock waves}
We can see that solution given by (\ref{quad1})-(\ref{quad2}) is, in general, multivalued. To figure out where the two-dimensional manifold $N$ given by (\ref{quad1})-(\ref{quad2}) has singularities of projection to the plane of independent variables, one needs to find zeroes of the two-form $dt\wedge dx$. Condition $(dt\wedge dx)|_{N}=0$ gives us a curve in the plane $\mathbb{R}^{2}(t,x)$ called \textit{caustic}. Choosing $\rho$ as a coordinate on the caustic, we get its equations in a parametric form:

\begin{equation}
\begin{split}
x(\rho)&=-\frac{1}{2\alpha_{1}}\left(2C_{1}\ln\rho+C_{1}(C_{3}^{3}\rho^{3}-4\rho^{2}C_{3}^{2}+3C_{3}\rho-2)+C_{3}C_{2}^{3}+2\alpha_{1}\alpha_{3}\right)\pm{}\\&\pm\frac{C_{2}(C_{3}\rho-1)^{2}}{\alpha_{1}\rho^{2}}\sqrt{C_{1}\rho^{3}+C_{5}\left(C_{3}+\frac{C_{7}}{\rho}\right)^{C_{6}}}-\frac{C_{5}\left(C_{3}+\frac{C_{7}}{\rho}\right)^{C_{6}}}{2(C_{6}+2)(C_{6}+3)C_{7}^3\alpha_{1}(C_{6}+1)\rho^{3}}\cdot{}\\&\cdot\left(C_{3}^{3}(-4+C_{7}^{3}(C_{6}^3+6C_{6}^2+11C_{6}+6)+(-2C_{6}-6)C_{7})\rho^3\right.-{}\\&-\left.2C_{7}((2(C_{6}^3+6C_{6}^2+11C_{6}+6))C_{7}^2+(-C_{6}^2-3C_{6})C_{7}-C_{6})C_{3}^2\rho^2+\right.{}\\&+\left.C_{7}^2(C_{6}+1)((C_{6}+3)(5C_{6}+12)C_{7}-2C_{6})C_{3}\rho-2C_{7}^3(C_{6}+4)(C_{6}+2)(C_{6}+1)\right),
\end{split}
\end{equation}

\begin{equation}
t(\rho)=-\alpha_{2}-\frac{C_{2}C_{3}}{\alpha_{1}}\pm\frac{(C_{3}\rho-1)^{2}}{\alpha_{1}\rho^{2}}\sqrt{C_{1}\rho^{3}+C_{5}\left(C_{3}+\frac{C_{7}}{\rho}\right)^{C_{6}}}.
\end{equation}

To construct a discontinuous solution from the multivalued one given by (\ref{quad1})-(\ref{quad2}), we use the mass conservation law. Let us write down equation (\ref{mass}) with the velocity $u$ found from (\ref{quad1}) in terms $t$ and $\rho$:

\begin{equation*}
\rho_{t}+\left(\rho\frac{\alpha_{1}\rho(t+\alpha_{2})+C_{2}}{1-C_{3}\rho}\right)_{x}=0,
\end{equation*}

and therefore the conservation law has the form

\begin{equation*}
\Theta=\rho dx-\rho\frac{\alpha_{1}\rho(t+\alpha_{2})+C_{2}}{1-C_{3}\rho}dt.
\end{equation*}

Its restriction $\Theta|_{N}$ to the manifold $N$ given by (\ref{quad1})-(\ref{quad2}) is a closed form, locally $\Theta|_{N}=dH$, and the potential $H(\rho,t)$ equals

\begin{equation*}
\begin{split}
H(\rho,t)&=\frac{\rho}{2\alpha_{1}(C_{3}\rho-1)^{2}}\left(C_{1}C_{3}^{3}\rho^{3}-4C_{1}C_{3}^{2}\rho^{2}+\rho\left( C_{2}^2C_{3}^2+(2C_{2}(t+\alpha_{2})\alpha_{1}+5C_{1})C_{3}+\alpha_{1}^{2}(t+\alpha_{2})^{2}\right)-2C_{1}\right)-{}\\&-\frac{C_{5}\left(C_{3}+\frac{C_{7}}{\rho}\right)^{C_{6}}}{(C_{6}+2)\alpha_{1}C_{7}^{2}(C_{6}+1)\rho^2}(C_{3}\rho+C_{7})(C_{3}(1+(C_{6}+2)C_{7})\rho-(C_{6}+1)C_{7}).
\end{split}
\end{equation*}

The discontinuity line, or a shock wave front is found from the system of equations

\begin{equation*}
H(\rho_{1},t)=H(\rho_{2},t),\quad x(\rho_{1},t)=x(\rho_{2},t),
\end{equation*}

where $x(\rho,t)$ is obtained from (\ref{quad1})-(\ref{quad2}) by eliminating $u$. Caustics along with the shock wave front are shown in Figure 3. Note that the picture is similar to that in case of phase transitions.

\begin{figure}[h!]
\centering
\includegraphics[scale=0.4]{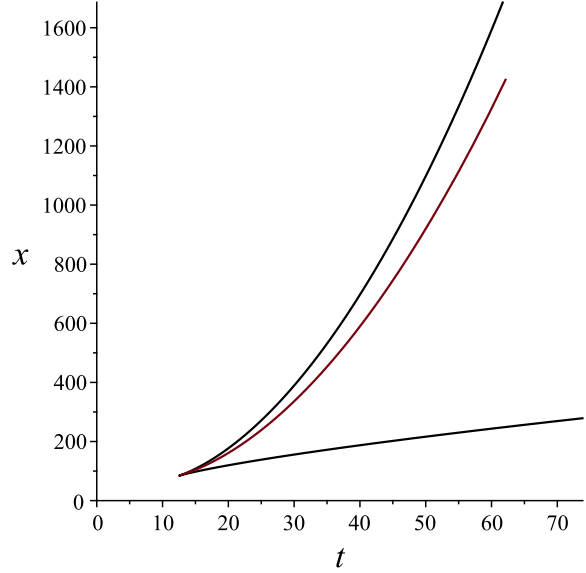}
\caption{Caustic (black) and shock wave front (red) for $n=3$.}
\end{figure}

The final result here is the expression for the time interval, within which the solution (\ref{quad1})-(\ref{quad2}) is smooth.
\begin{theorem}
Solution given by (\ref{quad1})-(\ref{quad2}) is smooth and unique in the time interval $t\in[0,t^{*})$, where

\begin{equation*}
t^{*}=\frac{1}{\alpha_{1}}\left(-C_{2}C_{3}-\alpha_{1}\alpha_{2}+(C_{3}-3)^{2}\sqrt{\frac{C_{1}}{27}+C_{5}(C_{3}+3C_{7})^{C_{6}}}\right),
\end{equation*}

and in case of (\ref{vdW-const}), where $C_{5}=240$, $n=3$, together with $C_{2}=1$, $\alpha_{1}=1$, $\alpha_{2}=2$, $\alpha_{3}=1$ approximately $t^{*}=12.53$.
\end{theorem}

\subsection{Phase transitions}
Having a solution, one can remove the phase transition curve from the space of thermodynamic variables to $\mathbb{R}^{2}(t,x)$. Indeed, on the one hand, we have all the thermodynamic parameters as functions of $(t,x)$. On the other hand, we have conditions on phase transitions (\ref{ph-tr-cond}) in the space of thermodynamic variables. In combination they give us a curve of phase transitions in $(t,x)$ plane. Phase transitions together with the shock wave are presented in Figure 4.

\begin{figure}[h!]
\centering
\includegraphics[scale=1]{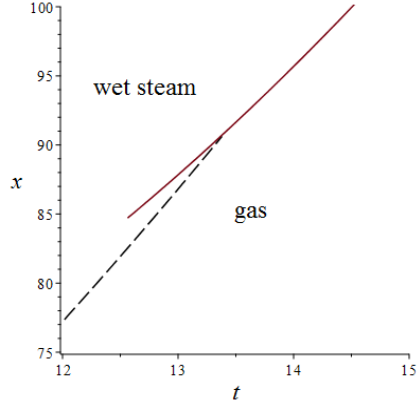}
\caption{Phase transition curve (dash line) and shock wave front (red line).}
\end{figure}

\section{Discussion}
In the present work we analyzed critical phenomena in gas flows of purely thermodynamic nature, which are phase transitions, and shock waves arising from singularities of solutions to the Euler system. To obtain such solutions we used a differential constraint compatible with the original system. In this work it was found in a purely computational way, and it seems interesting how to get it in a more regular way. One of possible ways to find such constraints is using differential invariants. Then, constraints can be found constructively by solving quotient PDEs, which was successfully realized in \cite{Eiv}. We hope to make use of this method in the future research. The analysis of phase transitions showed that sometimes shock waves can be accompanied with phase transitions, which is shown in Figure 4, since the phase transition curve intersects the shock wave front, and on the one side of the discontinuity curve we observe a pure gas phase, while on the other side we can see a wet steam.
\section*{Acknowledgements}
This work was partially supported by the Russian Foundation for Basic Research (project 18-29-10013) and by the Foundation for the Advancement of Theoretical Physics and Mathematics ``BASIS'' (project 19-7-1-13-3).


\begin{thebibliography}{99}
\bibitem{Arn1}
Arnold, V.  Singularities of Caustics and Wave Fronts. Springer Netherlands, 1990.
\bibitem{Arn2}
Arnold, V.  Catastrophe Theory. Springer-Verlag, Berlin, Heidelberg, 1984.
\bibitem{Arn3}
Arnold, V., Gusein-Zade, S., Varchenko, A.  Singularities of Differentiable Maps. Birkh\"{a}user, Basel, 1985.
\bibitem{Zeld}
Zeldovich, A., Kompaneets, I. Theory of Detonation. Academic Press, 1960.
\bibitem{Huang}
Huang, S.J., Wang, R. On blowup phenomena of solutions to the Euler equations for Chaplygin gases. {\em Applied Mathematics and Computation} {\bf 2013}, {\em 219}, 4365--4370.
\bibitem{Ros-Tab}
Rosales, R., Tabak, E. Caustics of weak shock waves. {\em Physics of Fluids} {\bf 1997}, {\em 10(1)}, 206--222.

\bibitem{Chat}
Chaturvedi, R., Gupta, P., Singh, L.P. Evolution of weak shock wave in two-dimensional steady supersonic flow in dusty gas. {\em Acta Astronautica} {\bf 2019}, {\em 160}, 552--557.
\bibitem{Polud1}
Poludnenko, A., Oran, E. The interaction of high-speed turbulence with flames: Global properties and internal flame structure. {\em Combustion and Flame} {\bf 2010}, {\em 157}, 995--1011.

\bibitem{Polud2}
Poludnenko, A., Gardiner, T., Oran, E. Spontaneous Transition of Turbulent Flames to Detonations in Unconfined Media. {\em Phys. Rev. Lett.} {\bf 2011}, {\em 107}, 054501.
\bibitem{LR-ljm-shock}
Lychagin, V., Roop, M. Shock waves in Euler flows of gases. {\em Lobachevskii J. Math.} {\bf 2020}, {\em 41(12)}, 2466-2472.
\bibitem{KLR}
Kushner, A., Lychagin, V., Rubtsov, V. Contact geometry and nonlinear differential equations. Cambridge University Press, Cambridge, 2007.

\bibitem{KrVin}
Vinogradov, A., Krasilshchik, I. (eds.) Symmetries and Conservation Laws for Differential Equations of Mathematical Physics. Factorial, Moscow, 1997.

\bibitem{KVL}
Vinogradov, A., Krasilshchik, I., Lychagin, V.  Geometry of jet spaces and nonlinear partial differential equations. Gordon and Breach, New York, 1996.

\bibitem{Ovs}
Ovsiannikov, L.  Group Analysis of Differential Equations. Academic Press, 1982.

\bibitem{Olver}
Olver, P.  Applications of Lie Groups to Differential Equations. Springer-Verlag, New York, 1986.
\bibitem{Tun}
Tunitsky, D. On multivalued solutions of equations of one-dimensional gas flow. {\em Proceedings of the 12th International Conference ``Management of Large-Scale System Development'' (MLSD)} {\bf 2019}, \url{https://ieeexplore.ieee.org/document/8911077}.

\bibitem{LychSing}
Lychagin, V. Singularities of multivalued solutions of nonlinear differential equations, and nonlinear phenomena. {\em Acta Appl.
Math.} {\bf 1985}, {\em 3(2)}, 135--173.
\bibitem{AKL-dan}
Akhmetzyanov, A., Kushner, A., Lychagin, V. Control of displacement front in a model of immiscible two-phase flow in porous media. {\em Doklady Mathematics} {\bf 2016}, {\em 94(1)}, 378--381.

\bibitem{AKL-ifac}
Akhmetzyanov, A., Kushner, A., Lychagin, V. Integrability of Buckley-Leverett's filtration model. {\em IFAC-PapersOnLine} {\bf 2016}, {\em 49(12)}, 1251--1254.

\bibitem{AKL-gsa}
Akhmetzyanov, A., Kushner, A., Lychagin, V. Shock waves in initial boundary value problem for filtration in two-phase 2-dimensional porous media. {\em Global and Stochastic Analysis} {\bf 2016}, {\em 3(2)}, 41--46.

\bibitem{KrugLych}
Kruglikov, B., Lychagin, V. Compatibility, Multi-Brackets and Integrability of Systems of PDEs. {\em Acta Appl. Math.} {\bf 2010}, {\em 109}, 151--196.
\bibitem{Eiv}
Schneider, E. Solutions of second-order PDEs with first-order quotients. \url{https://arxiv.org/abs/2005.06794}.
\bibitem{Lych-Yum-1}
Lychagin, V., Yumaguzhin, V. On Geometric Structures of 2-Dimensional Gas Dynamics Equations. {\em Lobachevskii J. Math.} {\bf 2009}, {\em 30(4)}, 327--332.

\bibitem{Lych-Yum-2}
Lychagin, V., Yumaguzhin, V. Minkowski Metrics on Solutions of the Khokhlov-Zabolotskaya Equation. {\em Lobachevskii J. Math.} {\bf 2009}, {\em 30(4)}, 333--336.
\bibitem{And}
Anderson, I., Torre, C.G. The differential geometry package. {\em Downloads.} {\bf 2016}, Paper 4. \url{http://digitalcommons.usu.edu/dg_downloads/4}.
\bibitem{Gibbs}
Gibbs, J.W. A Method of Geometrical Representation of the Thermodynamic Properties of Substances by Means of Surfaces {\em Transactions of the Connecticut Academy} {\bf 1873}, {\em 1}, 382--404.
\bibitem{Mrug}
Mrugala, R. Geometrical formulation of equilibrium phenomenological thermodynamics. {\em Reports on Mathematical Physics} {\bf 1978}, {\em 14(3)}, 419--427.
\bibitem{Rup}
Ruppeiner, G. Riemannian geometry in thermodynamic fluctuation theory. {\em Reviews of Modern Physics} {\bf 1995}, {\em 67(3)}, 605--659.

\bibitem{Lych}
Lychagin, V. Contact Geometry, Measurement, and Thermodynamics. In {\em Nonlinear PDEs, Their Geometry and Applications}; Kycia, R., Schneider, E., Ulan, M., Eds.; Birkh\"{a}user: Cham, Switzerland, 2019; pp. 3--52.
\bibitem{LR-ljm}
Lychagin, V., Roop, M. Critical phenomena in filtration processes of real gases. {\em Lobachevskii J. Math.} {\bf 2020}, {\em 41(3)}, 382-399.
\bibitem{KLR-ent}
Kushner, A., Lychagin, V., Roop, M. Optimal Thermodynamic Processes for Gases. {\em Entropy} {\bf 2020}, {\em 22(4)}, 448.
\end{thebibliography}
\end{document}